\documentstyle[twoside,12pt,toc_small,bib_toc,psfig]{article}

\newfont{\bbf}{msbm10 at 12pt}

\def\C{\mbox{\bbf C}}
\def\Cbar{\overline{\cz}}

\def\disk{\mbox{\,\bbf D}}

\def\diskbar{\ovl{\disk}}

\def\Circle{\mbox{\bbf S}^1}

\def\cz{\mbox{\bbf C}}

\def\ovl{\overline}
\def\phi1{\phi}
\def\phi{\varphi}
\def\eps{\varepsilon}

\def\*{ {\mbox{\tt $\star$}} }
\def\0{ {\mbox{\tt 0}} }
\def\1{ {\mbox{\tt 1}} }
\def\2{ {\mbox{\tt 2}} }
\def\3{ {\mbox{\tt 3}} }
\def\4{ {\mbox{\tt 4}} }

\def\St#1#2 {{\small\tt\rule{0pt}{0pt}_{\mbox{#1}}^{\mbox{#2}} }}
\def\<{\prec}
\def\>{\succ}

\ifx\ThmNum\undefined 
   \newtheorem{theorem}{Theorem}[section]
   \def\newsection#1{\setcounter{theorem}{0} \section{#1}}
\else 
   \newtheorem{theorem}{Theorem} 
   \def\newsection#1{\section{#1}}
\fi

\newtheorem{proposition}[theorem]{Proposition}
\newtheorem{lemma}[theorem]{Lemma}
\newtheorem{definition}[theorem]{Definition}

\newtheorem{corollary}[theorem]{Corollary}

\def\proof{\par\medskip\noindent {\sc Proof. }}
\def\proofof #1 {\par\medskip\noindent {\sc Proof of #1. }}
\def\sketch{\par\medskip\noindent{\sc Sketch of proof. }}
\def\sketchof #1 {\par\medskip\noindent {\sc Sketch of proof of #1. }}
\def\Box{\framebox[10pt]{\rule{0pt}{3pt}}}
\def\qed{\hfill $\Box$ \medskip \par}
\def\qedd{\rule{0pt}{0pt}\hfill $\Box$}
\def\remark{\par\medskip \noindent {\sc Remark. }}
\def\lineclear{\rule{0pt}{0pt}\par\noindent}

\def\LabelCaption#1#2{
   \centerline{\parbox{\captionwidth}{   %\vspace{-10mm} 
   \caption{\sl #2}  \label{#1} } }
}

\def\reminder #1 {{\sf #1}}
\def\hide #1 {}
\long\def\longhide #1 {}

\newfont{\script}{eusm10 at 12pt}
\def\M{{\mbox{\script M}}}  %    {{\bf M}}

\def\LabelCaption#1#2{
   \label{#1}
   \centerline{\parbox{\captionwidth}{   %\vspace{-10mm} 
   \caption{\sl #2}  } } }
\newfont{\bbflarge}{msbm10 at 21pt}

\def\A{{\cal A}}

\edef\theta{\vartheta}
\long\def\longhide#1{}

\def\MiLem{2.4}        % Mind. 3 Strahlen werden transitiv permutiert
\def\ExtRayLem{2.4}    % Das gleiche -- mit gleicher Nummer!
\def\LemOrbitSep{3.7}  % Orbit Separation Lemma
\def\LemFibers{2.4}    % Fundamental Properties of Fibers
\def\LemImpressFiber{2.5}  % Impression is in Single Fiber
\def\LemFibersNice{2.7}   % Nice fibers
\def\PropFiberLocConn{2.9} % Trivial Fibers imply loc conn
\def\PropLocConnFiber{2.10} % Loc conn makes fibers trivial
\def\ThmPeriodicFiber{3.5} % Repelling Periodic Points have Trivial F
\def\PropLocConnJuliaFiber{3.6} % 2.10 for Julias, specifies Q

\def\newsection#1{
   \section{#1}
}

\newlength\captionwidth
\title{\vspace{-15mm} On Fibers and Local Connectivity\\
of Mandelbrot and Multibrot Sets}
\author{Dierk Schleicher\\
Technische Universit\"at M\"unchen}
\date{}

\begin{document}

\maketitle
\thispagestyle{empty}
\def\IMSmarkvadjust{0 pt}
\def\IMSmarkhadjust{0 pt}
\def\SBIMSMark#1#2#3{
 \font\SBF=cmss10 at 10 true pt
 \font\SBI=cmssi10 at 10 true pt
 \setbox0=\hbox{\SBF Stony Brook IMS Preprint \##1}
 \setbox2=\hbox to \wd0{\hfil \SBI #2}
 \setbox4=\hbox to \wd0{\hfil \SBI #3}
 \setbox6=\hbox to \wd0{\hss
             \vbox{\hsize=\wd0 \parskip=0pt \baselineskip=10 true pt
                   \copy0 \break%
                   \copy2 \break% 
                   \copy4 \break}}
 \dimen0=\ht6   \advance\dimen0 by \vsize \advance\dimen0 by 8 true pt
                \advance\dimen0 by -\pagetotal
	        \advance\dimen0 by \IMSmarkvadjust
 \dimen2=\hsize \advance\dimen2 by .25 true in
	        \advance\dimen2 by \IMSmarkhadjust

%
%   Check for publication info
%
%  \newread\jref
  \openin2=publishd.tex
  \ifeof2\setbox0=\hbox to 0pt{}
  \else 
     \setbox0=\hbox to 3.1 true in{
                \vbox to \ht6{\hsize=3 true in \parskip=0pt  \noindent  
                {\SBI Published in modified form:}\hfil\break
                \input publishd.tex 
                \vfill}}
  \fi
  \closein2
  \ht0=0pt \dp0=0pt
 \ht6=0pt \dp6=0pt
 \setbox8=\vbox to \dimen0{\vfill \hbox to \dimen2{\copy0 \hss \copy6}}
 \ht8=0pt \dp8=0pt \wd8=0pt
 \copy8
 \message{*** Stony Brook IMS Preprint #1, #2. #3 ***}
}

\def\IMSmarkvadjust{-30pt}
\SBIMSMark{1998/13a}{December 1998}{}

%\centerline{\reminder{Draft of \today} }

\begin{center}
\begin{minipage}{110mm}
\def\toc_vspace{1.5em}
\def\tocname{Contents}
\tableofcontents{0.3em}
\end{minipage}
\end{center}
\vspace{5mm}

\begin{abstract} 
We give new proofs that the Mandelbrot set is locally connected at
every Misiurewicz point and at every point on the boundary of a
hyperbolic component. The idea is to show ``shrinking of puzzle
pieces'' without using specific puzzles. Instead, we introduce {\em
fibers} of the Mandelbrot set (see Definition~\ref{DefFiber}) and show
that fibers of certain points are ``trivial'', i.e., they consist of
single points.  This implies local connectivity at these points. 

Locally, triviality of fibers is strictly stronger than local
connectivity. Local connectivity proofs in holomorphic dynamics often
actually yield that fibers are trivial, and this extra knowledge is
sometimes useful. We include the proof that local connectivity of the
Mandelbrot set implies density of hyperbolicity in the space of
quadratic polynomials (Corollary~\ref{CorLocConnFiberHyp}).

We write our proofs more generally for {\em Multibrot sets}, which are
the loci of connected Julia sets for polynomials of the form $z\mapsto
z^d+c$. 

Although this paper is a continuation of \cite{Fibers}, it has been
written so as to be independent of the discussion of fibers of general
compact connected and full sets in $\C$ given there.
\end{abstract}

\newpage

\newsection {Introduction}
\label{SecIntro}

A great deal of work in holomorphic dynamics has been done in recent
years trying to prove local connectivity of Julia sets and of many
points of the Mandelbrot set, notably by Yoccoz, Lyubich, Levin,
van Strien, Petersen and others. One reason for this work is that the
topology of Julia sets and the Mandelbrot set is completely described
once local connectivity is known. Another reason is that local
connectivity of the Mandelbrot set implies that hyperbolicity is
dense in the space of quadratic polynomials.

In \cite{Fibers}, we have introduced {\em fibers} for arbitrary
compact connected and full subsets of the complex plane as a new
point of view on questions of local connectivity, and we have applied
fibers to filled-in Julia sets. In the present paper, we restrict to
the case of the Mandelbrot set and its higher degree cousins. After a
review of fundamental properties of the Multibrot sets, we define
fibers of these sets and show that they are born with much nicer
properties than fibers of general compact connected and full subsets
of $\C$. Our goal here is to prove that Misiurewicz points and
boundary points of hyperbolic components {\em have trivial fibers},
i.e.\ that their fibers consist of only one point. A consequence is
local connectivity at these points, but triviality of fibers is a
somewhat stronger property. Unlike in \cite{Fibers}, where the
discussion was for more general subsets of $\C$, we try to be as
specific as possible to Multibrot sets; we have repeated some key
discussions about fibers in order to make the present paper as
self-contained as possible.

We discuss polynomials of the form $z\mapsto z^d+c$, for
arbitrary complex constants $c$ and arbitrary degrees $d\geq 2$.
These are, up to normalization, exactly those polynomials which have
a single critical point. Following a suggestion of Milnor, we call
these polynomials {\em unicritical} (or unisingular). We will always
assume unicritical polynomials to be normalized as above, and the
variable $d$ will always denote the degree. We define the {\em
Multibrot set} of degree $d$ as the connectedness locus of these
families, that is
\[
\M_d := \{c\in\C: \mbox{ the Julia set of $z\mapsto z^d+c$ is
connected }\} \,\,.
\]
In the special case $d=2$, we obtain quadratic polynomials, and
$\M_2$ is the familiar {\em Mandelbrot set}. All the Multibrot sets
are connected, they are symmetric with respect to the real axis, and
they also have $d-1$-fold rotation symmetries (see
\cite{Intelligencer} with pictures of several of these sets). The
present paper can be read with the quadratic case in mind throughout.
However, we have chosen to do the discussion for all the Multibrot
sets because this requires only occasional slight modifications -- and
because recently interest in the higher degree case has increased; see
e.g.\ Levin and van Strien~\cite{LvS} or McMullen~\cite{McMullen}.

One problem in holomorphic dynamics is that many results are
folklore, with few accessible proofs published. In particular, many
of the fundamental results about the Mandelbrot set, due to Douady
and Hubbard, have been described in their famous ``Orsay notes''
\cite{Orsay} which have never been published. They are no longer
available, and it is not always easy to pinpoint a precise reference
even within these notes. This paper provides proofs of certain key
results about the Mandelbrot sets, and more generally for Multibrot
sets. Most of these results are known for $d=2$, but several of our
proofs are new and sometimes more direct that known proofs. 

In a sequel \cite{FiberTuning}, we will give various applications of
the concept of fibers which are related to the concepts of {\em
renormalization} and {\em tuning}: triviality of fibers is preserved
under tuning of Multibrot sets and under renormalization of Julia
sets, and we will show that the Mandelbrot set comes quite close to
being arcwise connected. These statements, in turn, have various
further applications. 

We begin in Section~\ref{SecMultibrot} by a review of certain
important properties of Mandelbrot and Multibrot sets. We then define
fibers for these sets in Section~\ref{SecFibers} using the same
definition as in \cite{Fibers}, except that the definition simplifies
because fibers of Multibrot sets generally behave quite nicely.

We also show that the fiber of an interior point is trivial if and
only if it is in a hyperbolic component. We conclude the section by a
proof that local connectivity is equivalent to triviality of all
fibers, and both conditions imply density of hyperbolicity (using an
argument of Douady and Hubbard; Corollary~\ref{CorLocConnFiberHyp}).

We then prove that every Multibrot set has trivial fibers and is thus
locally connected at every boundary point of a hyperbolic component
and at every Misiurewicz point. Boundaries of hyperbolic components
are discussed in Section~\ref{SecBoundary}, except roots of primitive
components: they require special treatment which can be found in
Section~\ref{SecRoots}. In Section~\ref{SecMisiurewicz}
we prove that fibers of Misiurewicz points are trivial. This shows
that fibers of Multibrot sets have particularly convenient properties.
Finally, in Section~\ref{SecCombinatorics}, we compare fibers to
combinatorial classes.

{\sc Acknowledgements.}
As a continuation of \cite{Fibers}, this paper inherits the
acknowledgements mentioned there. I would like to repeat the
contribution and support of Misha Lyubich in particular, as well as
the continuing encouragement of John Milnor and the support by the
Institute for Mathematical Sciences in Stony Brook.

\newsection{Multibrot Sets}
\label{SecMultibrot}

In this section, we review some necessary background about Multibrot
sets, and we define fibers of these sets. We include a result from
\cite{Fibers}: whenever a fiber consists of a single point, then
this implies local connectivity at this point.

All the Multibrot sets are known to be compact, connected and full. It
is conjectured but not yet known that they are locally connected.
However, it is known that many of its fibers are trivial. We will show
that for certain particularly important fibers in
Sections~\ref{SecBoundary}, \ref{SecRoots} and \ref{SecMisiurewicz}.

We recall the definition of external rays: For any compact connected
and full subset $K\subset\C$ consisting of more than a single point,
there is a unique conformal isomorphism
$\Phi\colon\Cbar-K\to\Cbar-\diskbar$ fixing $\infty$, normalized so as
to have positive real derivative at $\infty$. Inverse images of radial
lines in $\C-\diskbar$ are called external rays, and an external ray
at some angle $\theta$ is said to {\em land} if the limit
$\lim_{r\searrow 1} \Phi^{-1}(re^{i\theta})$ exists. The {\em
impression} of this external ray is the set of all limit points of
$\Phi^{-1}(r'e^{i\theta'})$ for $r'\searrow 1$ and $\theta'\to\theta$.

As in \cite{ExtRays} and \cite{MiOrbits}, we will denote external
rays of the Multibrot sets by {\em parameter rays}\/ in order to
distinguish them from {\em dynamic rays}\/ of Julia sets. All the
parameter rays at rational angles are known to land (see Douady and
Hubbard~\cite{Orsay}, Schleicher~\cite{ExtRays} or, in the periodic
case, Milnor~\cite{MiOrbits}). 

A {\em ray pair} is a collection of two external rays (dynamic or
parameter rays) which land at a common point. A dynamic ray pair is
{\em characteristic} if it separates the critical value from the
critical point and from all the other rays on the forward orbit of the
ray pair. The landing point of a periodic or preperiodic dynamic ray
pair is always on a repelling or parabolic orbit; if a preperiodic
dynamic ray pair is characteristic, then its landing point is
necessarily on a repelling orbit.

The following theorem will be used throughout this paper. The first
half is due to Lavaurs~\cite{La}. Proofs of both parts can be found in
\cite{MandelStruct}.

\begin{theorem}[Correspondence of Ray Pairs]
\label{ThmRayCorrespondence} \lineclear
For every degree $d\geq 2$ and every unicritical polynomial $z\mapsto
z^d+c$, there are bijections 
\begin{itemize}
\item  
between the ray pairs in parameter space at periodic angles,
separating $0$ and $c$, and the characteristic periodic ray pairs in
the dynamic plane of $c$ landing at repelling orbits; and
\item  
between the ray pairs in parameter space at preperiodic angles,
separating $0$ and $c$, and the characteristic preperiodic ray pairs
in the dynamic plane of $c$.
\end{itemize}  
These bijections of ray pairs preserve external angles.
\qedd
\end{theorem}
We assume that the separating ray pairs do not go through the point
$c$ (the critical value or the parameter). For these boundary cases,
we have the following. The critical value is never on a periodic ray
pair because, depending on whether the corresponding Julia set is
connected or not, it would either be periodic and thus superattracting
and could not be the landing point of dynamic rays, or the periodic
dynamic rays would bounce into a precritical point infinitely often
and thus fail to land. Conversely, if the parameter $c$ is on a
parameter ray pair at periodic angles, the corresponding Julia set
either has a dynamic ray pair which is the characteristic ray pair of
a parabolic periodic orbit, or at least one of the rays on this ray
pair bounces into a precritical point infinitely often. In the
preperiodic case, if the parameter $c$ is on a parameter ray pair at
preperiodic angles, then in the corresponding dynamic plane, the
dynamic rays at the corresponding angles form a dynamic ray pair
landing at a preperiodic point in a repelling orbit. Moreover, this
ray pair contains the critical value, and all its forward images are
on the same side as the critical point. Conversely, if in some dynamic
plane the critical value is on a characteristic preperiodic ray pair,
then this ray pair lands on a repelling orbit and the parameter is on
the parameter ray pair at corresponding angles. In the quadratic
case, these results go back to Douady and Hubbard~\cite{Orsay}; for
proofs, see also \cite[Theorem~1.1]{ExtRays} and, in the periodic case,
\cite{MiOrbits}. Proofs for higher degrees will be in \cite{Dominik}.

Not all rational parameter rays are organized in pairs. The number of
parameter rays at preperiodic angles landing at a common point can be
any positive integer. For parameter rays at periodic angles, this
number is either $1$ or $2$; in the quadratic case, this number is
always $2$ (we count the parameter rays at angles $0$ and $1$
separately).

All parameter rays at periodic angles are known to land at parabolic
parameters (where the critical orbit converges to a unique parabolic
orbit). All parameter rays at preperiodic angles land at parameters
where the critical value is strictly preperiodic: such parameters are
(somewhat unfortunately) known as ``Misiurewicz points''. Another name
for such points would be ``critically preperiodic parameters''.

Any hyperbolic component has one {\em root} and $d-2$ {\em
co-roots} on its boundary: these are parameters with parabolic orbits
whose periods are at most that of the attracting orbits within the
component. The root is the landing point of two periodic parameter
rays, while every co-root is the landing point of a single parameter
ray at a periodic angle. The period of a hyperbolic component is also
the period of the $d$ parameter rays landing at root and co-roots of
the component. 

The {\em wake} of a hyperbolic component will be the region in the
complex plane separated from the origin by the two rational parameter
rays landing at the root of the component. {\em Subwakes} of a component
will be wakes of components bifurcating directly from the component. The
intersection of any wake or subwake with the Multibrot set will be
called a {\em limb} or {\em sublimb}. Wakes and subwakes are open and
connected, while limbs and sublimbs are connected and closed, except
that one boundary point is missing (the ``root'' of the (sub-)limb,
where it is attached to the rest of the Multibrot set). Every hyperbolic
component has a unique {\em center:} a parameter where the critical
orbit is periodic. These centers and Misiurewicz points together form
the countable set of {\em postcritically finite} parameters. When
treating them simultaneously, it will sometimes simplify language to
identify the component with its center and speak, e.g.\ of the ``wake
of the center'', meaning of course the wake of the component containing
this center.

All the rational rays landing at any Misiurewicz point cut the
complex plane into as many open parts as there are rays. We will call
these parts the {\em subwakes} of the Misiurewicz point. The subwake
containing the origin will be called the {\em zero wake} (which is
the entire complex plane minus the ray if there is only one ray),
and the union of all the other wakes (together with any parameter
rays between them) will be called the {\em wake} of the Misiurewicz
point. 

Of fundamental importance for the investigation of the Multibrot sets
is the following Branch Theorem, which was first proved by Douady and
Hubbard for the Mandelbrot set: this theorem is one of the principal
results of their theory of ``nervures'' (Expos\'es XX--XXII in
\cite{Orsay}; the Branch theorem is their Proposition~XXII.3). Another
proof can be found in \cite[Theorem~9.1]{IntAddr}.

\begin{theorem}[Branch Theorem] 
\label{ThmBranch} \lineclear   
For every two postcritically finite parameters $c_1$ and $c_2$,
either one of them is in the wake of the other, or there is a
Misiurewicz point such that $c_1$ and $c_2$ are in two different of
its subwakes, or there is a hyperbolic component such that $c_1$ and
$c_2$ are in two different of its subwakes.
\end{theorem}

\begin{figure}[htbp]
%\vspace{50mm}
%\hspace{30mm}\special{pictfile P_ThreeRaysOneFiber}
\centerline{
\psfig{figure=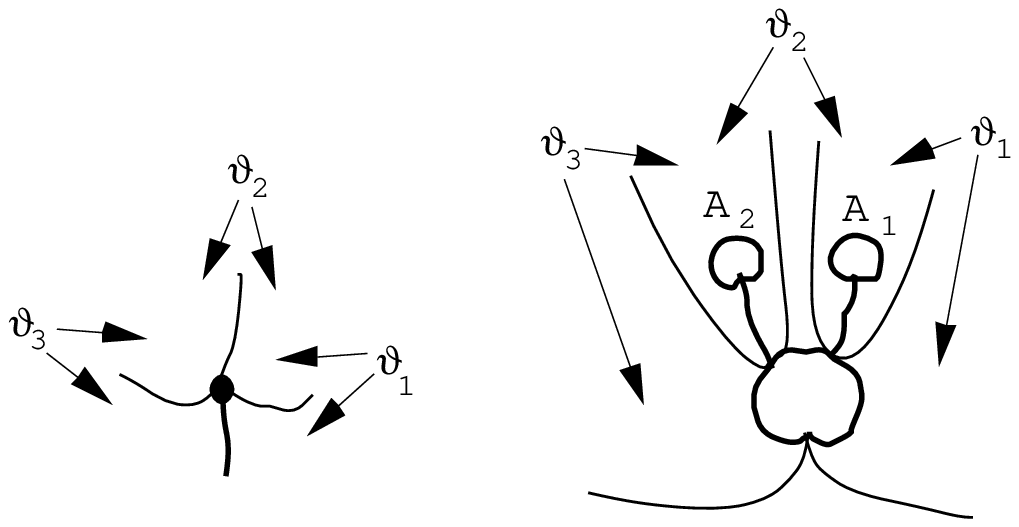,height=60mm}
}
\LabelCaption{FigThmBranch}
{Illustration of the Branch Theorem. 
Left: a Misiurewicz point with three parameter rays landing, and the
possible positions for the $\theta_i$; right: separation at a
hyperbolic component. }
\end{figure}

\proof    
We will work in the dynamical plane of $c_1$, which we tried to
sketch in Figure~\ref{FigMandelBranchDyn}. Let $\alpha$ be one of the
external angles of $c_2$ and let $a_2$ be the (pre-)periodic point at
which the dynamic ray at angle $\alpha$ lands. 

If the point $c_1$ is on the arc $[0,a_2]$, then two dynamic rays
land at $c_1$ (if it is preperiodic) or at the dynamic root of the
Fatou component containing $c_1$ (if it is periodic) and, by
Theorem~\ref{ThmRayCorrespondence}, the parameter rays at the same
angles land at $c_1$ (resp.\ the root of the hyperbolic component
containing $c_1$) and separate the parameter ray at angle $\alpha$
and the point $c_2$ from the origin, so $c_2$ is in the wake of $c_1$.
Similarly, if $a_2$ is on the arc $[0,c_1]$, then $c_1$ is in the wake
of $c_2$. If none of these cases occurs, then the intersection of the
arcs $[0,a_2]$ and $[0,c_1]$ is the arc $[0,b]$ for some point $b$
which is a branch point in the Julia set. The union of forward images
of the regular arcs $[0,c_1]$ and $[0,a_2]$ has the topological type
of a finite tree; in particular, it has finitely many branch points.
Since this tree is forward invariant, all these branch points are
periodic or preperiodic (or critical, but then they have finite
forward orbits by assumption); in particular, the forward orbit of
$b$ is finite.

\begin{figure}[htbp]
%\vspace{70mm}
%\hspace{65mm}\special{pictfile MandelBranchDyn.pict scaled 400}
\centerline{
\psfig{figure=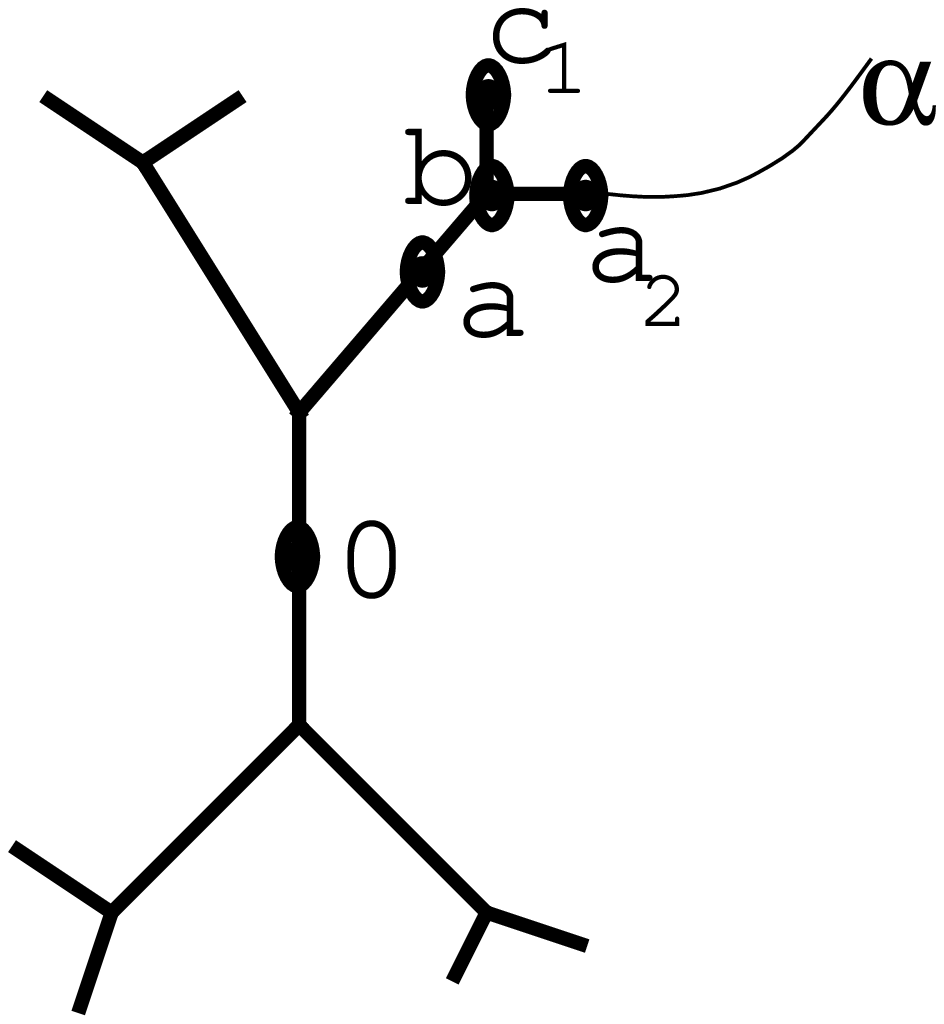,height=50mm}
}
\LabelCaption{FigMandelBranchDyn}
{A sketch of the dynamic plane of the parameter $c_1$ as used in the
proof of the Branch Theorem~\protect\ref{ThmBranch}.} 
\end{figure}

Now consider the set $Z$ of characteristic (pre-)periodic points on
$[0,b]$ and let $a$ be the supremum of $Z-\{b\}$ on the arc $[0,b]$
with respect to the induced order on this arc. It is slightly easier
to consider the case that $c_1$ is a Misiurewicz point, so every
periodic and preperiodic point is on a repelling orbit. If $c_1$ is
a center of a hyperbolic component, the point $a$ cannot be on the
superattracting periodic orbit because the Fatou component it is in
separates the critical value from the origin, so $a$ cannot be a
characteristic point. Therefore, $a$ is on a repelling orbit and thus
the landing point of at least two dynamic rays. For a similar reason,
if the point $b$ is characteristic, it cannot be on the
superattracting orbit.

%\begin{description}
%\item[

\medskip

{\em The case that $a$ and $b$ are two different points:}\/
%] \lineclear
if the arc $[a,b]$ contains periodic points on its interior, let $z$
be one of lowest period; if there is no periodic point on the
interior of this arc, let $z:=b$. Then the point $z$ is periodic or
preperiodic; assuming for the moment that it not characteristic,
then there is a first number $s$ of iterations after which the point
$z$ maps behind itself. If $z$ is periodic, then $s$ must be smaller
than the period. Under $s$ iterations, the image of the arc $[a,z]$
contains the arc itself. By the intermediate value theorem, there is a
point on $[a,z]$ which maps to itself after $s$ iterations. This
point cannot be $z$ by construction, and neither can it be an interior
point of the arc. It therefore follows that the point $a$ is periodic,
and its period is $s$ or divides $s$. It is the landing point of at
least two dynamic rays, and the two characteristic rays reappear as
parameter rays landing at the root of a hyperbolic component by
Theorem~\ref{ThmRayCorrespondence}. Looking at external angles, it
follows that the points $c_1$ and $c_2$ are contained in the wake of
this component. If they are contained within the same subwake, the two
parameter rays bounding this subwake would correspond to a periodic
point in $Z$ behind $a$, contradicting maximality of $a$. 

If, however, $z$ is characteristic, then $z=b$ and there is no
periodic point on the arc $[a,b]$. The construction above goes
through except if $a$ and $b$ are periodic and, after a finite number
$s$ of iterations, the arc $[a,b]$ covers itself homeomorphically.
But then $a$ and $b$ cannot be both repelling as they must be.

\medskip

%\item[
{\em The case that $a$ and $b$ are equal:}\/
%] \lineclear
in this case, at least three dynamic rays land at $a=b$, and
these rays separate the points $c_1$ and $a_2$ from each other and
from the origin. The point $b$ is characteristic: otherwise, it would
map behind itself after finitely many steps; by continuity, it could
not be the limit of characteristic points in $Z$. 

The point $b$ is a limit of points in $Z$. If $b$ is periodic, the
dynamic rays landing at $b$ would be permuted transitively by the
first return map of $b$, and points in $Z$ sufficiently close to
$b$ would be mapped behind $b$ and thus behind themselves. This is a
contradiction. Therefore, $b$ is preperiodic, and 
Theorem~\ref{ThmRayCorrespondence} turns the three dynamic rays
landing at $b$ into three parameter rays landing at a common
Misiurewicz point such that $c_1$ and $c_2$ are in two different of
its subwakes because these subwakes contain the parameter rays
associated to $c_1$ and $c_2$. 
%\end{description}
\medskip

This finishes the proof of the theorem.
\qed
%\reminder{Adapt the proof to current terminology!}

\remark
In \cite{FiberTuning}, we show that the Mandelbrot set is ``almost''
arcwise connected. The Branch Theorem implies then that arbitrary
``regular'' arcs in the Mandelbrot set connecting postcritically
finite parameters can branch off from each other only at Misiurewicz
points or within hyperbolic components.

As a first corollary, we can describe how many rational parameter rays
may land at the boundary of any interior component of a Multibrot set.
All the known interior components of Multibrot sets are ``hyperbolic''
components in which the dynamics has an attracting periodic orbit. It
is conjectured that non-hyperbolic (``queer'') components do not exist.

\begin{corollary}[Interior Components of the Multibrot Sets]
\label{CorMandelInterior} \lineclear
Every hyperbolic component of a Mandelbrot or Multibrot set has
infinitely many boundary points which are landing points of rational
parameter rays, and all of them are accessible from inside the
component. However, every non-hyperbolic component has at most
one boundary point which is the landing point of a parameter ray at a
rational angle
\end{corollary}
\remark
In Section~\ref{SecCombinatorics}, we will strengthen the second
statement by showing that no rational parameter ray can ever land on
the boundary of a non-hyperbolic component (provided such a thing
exists at all). 
\proof
The statement about hyperbolic components is well known at least in
the quadratic case; see Douady and Hubbard~\cite{Orsay},
Milnor~\cite[Section~7]{MiOrbits},  or
Schleicher~\cite[Section~5]{ExtRays}. For the general case, see
\cite{Dominik}.

Assume that two rational parameter rays land at different points on
the boundary of a non-hyperbolic component. Denote their external
angles by $\theta_1$ and $\theta_2$ and the landing points by $c_1$
and $c_2$. For readers familiar with kneading sequences and their
geometric interpretation as internal addresses \cite{IntAddr}, there
is a brief argument: the external angles $\theta_i$ translate into
kneading sequences and angled internal addresses. If the latter are 
different, then the landing points of the two parameter rays must be
separated by a parameter ray pair at periodic angles (this argument,
without using angles in the internal address, was used in Sections~3
and 4 of \cite{ExtRays} to establish the landing properties of
parameter rays we are using here). However, if the two angled
internal addresses coincide, then the two parameter rays land at a
common point \cite[Theorem~9.2]{IntAddr}.

Here we give a different proof, not involving internal addresses.
If one of the $\theta_i$ is periodic, then replace the corresponding
parabolic point $c_i$ by the center of the component at the root or
co-root of which the ray at angle $\theta_i$ lands. Since $c_1$ and
$c_2$ are on the closure of the same non-hyperbolic component, there
cannot be a Misiurewicz point or a hyperbolic component separating
$c_1$ and $c_2$ from each other and from the origin. Therefore, by the
Branch Theorem~\ref{ThmBranch}, one of these two points must be within
the wake associated to the other point. Without loss of generality,
assume that  $c_2$ is within the wake of the Misiurewicz point $c_1$ or
of the hyperbolic component with center $c_1$. Take a third preperiodic
angle $\theta_3$ within the wake of $c_1$ but outside of the wake of
$c_2$, and so that its landing point $c_3$ is a Misiurewicz point
different from $c_1$ and $c_2$. If $c_2$ is within the wake of the
Misiurewicz point $c_3$, then two rays landing at this Misiurewicz point
separate $c_1$ and $c_2$, so these two points cannot be on the boundary
of the same non-hyperbolic component. Otherwise the Branch Theorem
supplies another Misiurewicz point or hyperbolic component separating
$c_2$ and $c_3$ from each other and from the origin. If the angle
$\theta_3$ is chosen sufficiently closely to $\theta_2$, this
separation point cannot be $c_1$, so this new Misiurewicz point or
hyperbolic component separates $c_2$ from $c_1$, and again these two
points cannot be on the boundary of a common non-hyperbolic component. 
\qed

The most important application of the Branch Theorem is to relate
local connectivity and density of hyperbolicity for Mandelbrot and
Multibrot sets. We will do that in Corollary~\ref{CorLocConnFiberHyp},
showing also that local connectivity and triviality of all fibers is
equivalent.

\newsection{Fibers of Multibrot Sets}
\label{SecFibers}

We will now define fibers of Multibrot sets using parameter rays
at rational angles; recall that all these rays are known to land. We
first need to introduce separation lines.

\begin{definition}[Separation Line]
\label{DefSeparation} \lineclear
A {\em separation line} is either a pair of parameter rays at rational
angles landing at a common point (a ray pair), or a pair of parameter
rays at rational angles landing at the boundary of the same interior
component of the Multibrot set, together with a simple curve within
this interior component which connects the landing points of the two
rays. Two points $c,c'$ in a Multibrot {\em can be separated} if there
is a separation line $\gamma$ avoiding $c$ and $c'$ such that\/ $c$
and\/ $c'$ are in different connected components of\/ $\C-\gamma$.
\end{definition}

Of course, the interior components used in this definition will always
be hyperbolic components by Corollary~\ref{CorMandelInterior}.

It will turn out that points on the closure of a hyperbolic component
of a Multibrot set have trivial fibers. For all other points, it will
be good enough to construct fibers using ray pairs at rational
external angles (in fact, periodic external angles suffice). We will
justify this in Proposition~\ref{PropPeriodicParaRays}.

\begin{definition}[Fibers and Triviality]
\label{DefFiber} \lineclear
The fiber of a point\/ $c$ in a Multibrot set\/ $\M_d$ is the set of
all points $c'\in\M_d$ which cannot be separated from $c$. The fiber
of $c$ is {\em trivial} if it consists of the point $c$ alone.
\end{definition}

A fundamental construction in holomorphic dynamics is called the {\em
puzzle}, introduced by Branner, Hubbard and Yoccoz. A typical proof
of local connectivity consists in establishing {\em shrinking of
puzzle pieces} around certain points. This is exactly the model for
fibers: the fiber of a point is the collection of all points which
will always be in the same puzzle piece, no matter how the puzzle was
constructed. Our arguments will thus never use specific puzzles. 
The definition here is somewhat simpler than in \cite{Fibers}. The
reason is that the definition in \cite{Fibers} is for arbitrary
compact connected and full subsets of $\C$, while some conceivable
difficulties cannot occur for Multibrot sets, due to the Branch
Theorem~\ref{ThmBranch} and its Corollary~\ref{CorMandelInterior}.
(The equivalence of both definitions follows from this Corollary and
\cite[Lemma~\LemFibers]{Fibers}.)

\begin{proposition}[Fibers of Interior Components]
\label{PropFibersInt} \lineclear
The fiber of a point within any hyperbolic component of a Multibrot
set is always trivial. The fiber of a point within any non-hyperbolic
component of a Multibrot set always contains the closure of its
non-hyperbolic component.
\end{proposition}
\proof
Let $c$ be any interior point in a hyperbolic component. It can easily
be separated from any other point within the same hyperbolic component
by a separation line consisting of two rational parameter rays landing
on the boundary of this hyperbolic components, and a curve within this
component; there are infinitely many such rays by
Corollary~\ref{CorMandelInterior}. 

Any boundary point of the hyperbolic component, and any point in
$\M_d$ outside of this component, can also be separated from $c$ by
such a separation line (here it is important to have more than two
boundary points of the component which are landing points of rational
parameter rays; if there were only two such boundary points, then
these points could not be separated from any interior point). 

Since any non-hyperbolic component has at most one boundary point
which is the landing point of rational rays, no separation line can
run through such a component, and the fibers of all such interior
points contain at least the closure of this non-hyperbolic component.
\qed

Here are a few more useful properties of fibers. They are all taken
from \cite[Lemma~\LemFibers]{Fibers}, using
Corollary~\ref{CorMandelInterior} to exclude certain bad
possibilities. For convenience, we repeat most proofs in our context.

\begin{lemma}[Properties of Fibers]
\label{LemFibers} \lineclear
Fibers have the following properties:
\begin{enumerate}
\item
Fibers are always compact, connected and full. 
\item
The relation ``is in the fiber of'' is symmetric: for two points
$c,c'\in\M_d$, either each of them is in the fiber of the other one,
or both points can be separated.
\item
The boundary of any fiber is contained in the boundary of\/ $\M_d$,
unless the fiber is trivial.
\end{enumerate}
\end{lemma}
\proof
If $c'$ is not in the fiber of $c$, then by definition these two
points can be separated, and $c$ is not in the fiber of $c'$. This
proves the second claim.

No interior point of $\M_d$ can be in the boundary of any non-trivial
fiber: a hyperbolic interior point is a fiber by itself, and closures
of non-hyperbolic components are completely contained in fibers. This
settles the third claim.

Every fiber is compact because every separation line separates an open
subset of $\M_d$. In the construction of the fiber of any point
$c\in\M_d$, ray pairs always leave connected and full neighborhoods of
$c$, and nested intersections of compact connected and full subsets of
$\C$ are compact connected and full. The only further possible
separation lines run through hyperbolic components, but hyperbolic
interior points have trivial fibers, and we only have to consider 
boundary points of hyperbolic components. These points can be treated
easily; details can be found in \cite[Lemma~\LemFibers]{Fibers}.
(Equivalently, we will prove in Sections~\ref{SecBoundary} and
\ref{SecRoots} that boundary points of hyperbolic components have
trivial fibers.)
\qed

A key observation is that triviality of a
fiber implies that $\M_d$ is locally connected at this point. We
repeat the argument from \cite[Proposition~\PropFiberLocConn]{Fibers}.

\noindent
\begin{minipage}{\textwidth}
\begin{proposition}[Trivial Fibers Yield Local Connectivity]
\label{PropFiberLocConn} \lineclear  
If the fiber of a point $c\in\M_d$ is trivial, then $\M_d$ is openly
locally connected at $c$. Moreover, if the parameter ray at some angle
$\theta$ lands at a point $c$ with trivial fiber, then for any
sequence of external angles converging to $\theta$, the corresponding
impressions converge to $\{c\}$. In particular, if all the fibers of
$\M_d$ are trivial, then $\M_d$ is locally connected, all external
rays land, all impressions are points, and the landing points depend
continuously on the angle.
\end{proposition}
\end{minipage}
\medskip

\proof 
Consider a point $c\in\M_d$ with trivial fiber. If $c$ is in the
interior of $\M_d$, then $K$ is trivially openly locally connected at
$c$. Otherwise, let $U$ be an open neighborhood of $c$. By
assumption, any point $c'$ in $\M_d-U$ can be separated from $c$ such
that the separation avoids $c$ and $c'$. The region cut off from $c$
is open; what is left is a neighborhood of $c$ having connected
intersection with $\M_d$. By compactness of $\M_d-U$, a finite number
of such cuts suffices to remove every point outside $U$, leaving
another neighborhood of $c$ intersecting $\M_d$ in a connected set.
Removing the finitely many cut boundaries, an open neighborhood
remains, and $\M_d$ is openly locally connected at $c$. Similarly, if
$c$ is the landing point of the $\theta$-ray, then external rays with
angles sufficiently close to $\theta$ will have their entire
impressions in $\ovl U$ (although the rays need not land). 
Finally, it is easy to see that the impression of any ray is contained
in a single fiber (compare \cite[Lemma~\LemImpressFiber]{Fibers}).
\qed

Now, we show that the fiber of an interior point is trivial if and
only if it is within a hyperbolic component, and that local
connectivity of a Multibrot set is equivalent to all its fibers
being trivial. Both properties will imply that every interior
component is hyperbolic.

We use the ideas of the original proof of Douady and
Hubbard~\cite[Expos\'e~XXII.4]{Orsay}.

\begin{corollary}[Local Connectivity, Trivial Fibers and Hyperbolicity] 
\label{CorLocConnFiberHyp} \lineclear
For the Mandelbrot and Multibrot sets, local connectivity is
equivalent to triviality of all fibers, and both imply density of
hyperbolicity.
\end{corollary}

\remark
Local connectivity is definitely stronger than density of
hyperbolicity: the former amounts to fibers being points, while the
latter means only that fibers have no interior. This argument is taken
from Douady~\cite{DoCompacts}, where it is part of a sketch of an
independent proof. In our context, it is just a restatement of the
original proof in \cite{Orsay}.

\proof
We know from Proposition~\ref{PropFiberLocConn} that triviality of all
fibers implies local connectivity. Moreover, any non-hyperbolic
component would be contained in a single fiber, so triviality of all
fibers implies that all interior components are hyperbolic. Since the
exterior of $\M_d$ is hyperbolic as well, triviality of all fibers
of $\M_d$ implies that hyperbolic dynamics is dense in the space of
all unicritical polynomials of degree $d$. 

It remains to show that local connectivity of a Multibrot set implies
that all fibers are trivial. For this, we assume that the Multibrot
set $\M_d$ is locally connected. If there is a fiber which is not a
singleton, or which even contains a non-hyperbolic component, denote
it $Y$. The set $Y$ is then uncountable and its boundary is contained
in the boundary of $\M_d$ by Lemma~\ref{LemFibers}. Let $c_1,c_2,c_3$
be three boundary points which are not landing points of any of the
countably many parameter rays at rational angles. By local
connectivity and Carath\'eodory's Theorem, there are three parameter
rays at angles $\theta_1,\theta_2,\theta_3$ landing at these points,
and these three rays separate $\C-Y$ into three regions. Similarly,
the three angles cut $\Circle$ into three open intervals. Pick one
periodic angle from each of the  two intervals which do not contain
$0$; the corresponding parameter rays land at roots or co-roots of two
hyperbolic components $\A_1$ and $\A_2$. The three parameter rays at
angles $\theta_i$, together with $Y$, separate these two components
from each other and from the parameter ray at angle $0$ and thus from
the origin. Applying the Branch Theorem~\ref{ThmBranch} to these two
components, there must either be a Misiurewicz point or a hyperbolic
component separating these two components from each other or from the
origin. If it is a Misiurewicz point, then three rational rays landing
at it must separate the three points $c_1,c_2,c_3$ from each other,
which is incompatible with $Y$ being a single fiber or a
non-hyperbolic component. Similarly, if the separation is given by a
hyperbolic component, then this component, together with the parameter
rays forming its wake and the subwakes containing $\A_1$ and $\A_2$,
again separate the three $c_i$, yielding the same contradiction.

We conclude that, if a Multibrot set is locally connected, then
its fibers are trivial, and every connected component of its interior
is hyperbolic.
\qed

\begin{figure}
%\vspace{50mm}
%\hskip 50mm
%\special{pictfile LocConnFiber.pict scaled 500}
\centerline{
\psfig{figure=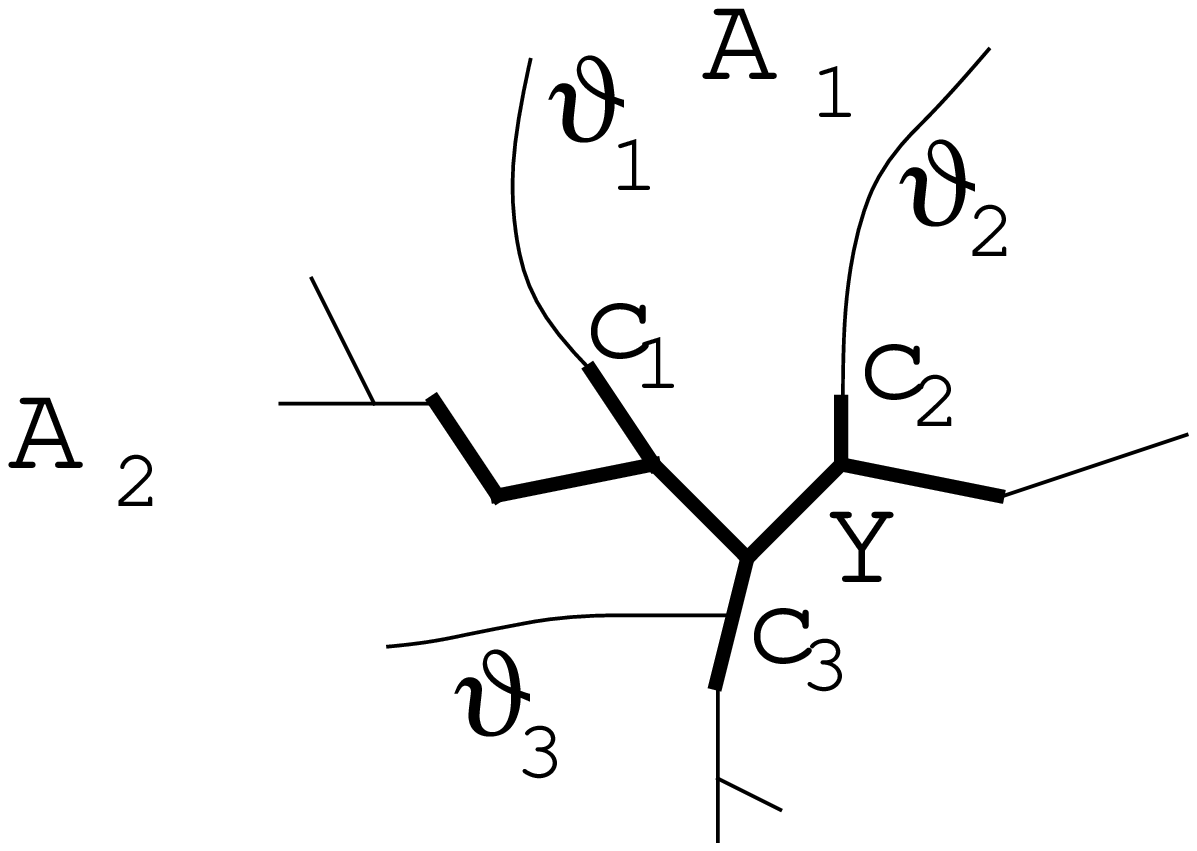,width=80mm}
}
\LabelCaption{FigMLC_HD}
{Illustration of the proof that local connectivity implies that all
fibers are trivial. Highlighted is a non-trivial fiber $Y$, and three
of its boundary points $c_1,c_2,c_3$ together with parameter rays
landing there. The Branch Theorem, applied to the hyperbolic
components $\A_1$ and $\A_2$, provides then a contradiction.
}
\end{figure}

\remark
In \cite[Proposition~\PropLocConnFiber]{Fibers}, we have shown that
local connectivity of any compact connected and full subset of $\C$ is
equivalent to triviality of all fibers, provided that fibers are
constructed using an appropriate collection of external rays. We have
given an extra proof here in order to make this paper self-contained,
and in order to show that the external rays at rational angles
suffice. We will show in Section~\ref{SecCombinatorics} that even
parameter rays at periodic external angles are sufficient.

\newsection{Boundaries of Hyperbolic Components}
\label{SecBoundary}

In this section, we will study boundary points of hyperbolic
components, except roots of primitive components which will require
special treatment (see Section~\ref{SecRoots}). The following result
was originally proved by Yoccoz using a bound on sizes of sublimbs
(the ``Yoccoz inequality''); see Hubbard~\cite[Section~I.4]{HY}. That
proof, as well as ours, uses in an essential way the fact that
repelling periodic points in connected Julia sets are landing points
of periodic dynamic rays. This fact is due to Douady or Yoccoz
\cite[Theorem~I.A]{HY}.

\begin{theorem}[No Irrational Decorations] \label{ThmNoIrratDecos}
\lineclear  
Any point in $\M_d$ within the wake of a hyperbolic component is
either in the closure of the component or within one of its sublimbs
at rational internal angles with denominator at least two.
\end{theorem}
\remark
Sometimes, this theorem is phrased as saying that hyperbolic
components have ``no ghost limbs'': decorations are attached to
hyperbolic components only at parabolic boundary points, but not at
root, co-roots, or irrational boundary points.

\proof 
Let $W$ be a hyperbolic component, let $n$ be its period and let
$\tilde c$ be a point within the limb of the component but not on the
closure of $W$. There are finitely many hyperbolic components of
periods up to $n$, some of which might possibly be within the wake of
$W$. We lose nothing if we assume $\tilde c$ to be outside the closures
of the wakes of such hyperbolic components, possibly after replacing
$\tilde c$ with a different parameter: if there was a hyperbolic
component of period up to $n$ in an ``irrational sublimb'' of $W$, then
this same ``irrational sublimb'' is connected and thus contains points
arbitrarily close to $W$. More precisely, we can argue as follows: the
rational parameter rays landing on the boundary of $W$ split into two
groups according to whether they pass $\tilde c$ to the ``left'' or
``right'', and the region in the wake of $W$ sandwiched between these
rays is connected. (Moreover, it is quite easy to show that every
hyperbolic component in the wake of $W$ is contained in a subwake at
rational internal angle; see \cite{MandelStruct}). 

Denote the center of $W$ by $c_0$. There is a path $\gamma$
connecting $c_0$ to $\tilde c$ within the wake of $W$ and avoiding
the closures of the wakes of all hyperbolic components of periods up to
$n$ within the wake of $W$, except $W$ itself. This path need not be
contained within the Multibrot set. For the parameter $c_0$, the
critical value is a periodic point. We want to continue this periodic
point analytically along $\gamma$, obtaining an analytic function
$z(c)$ such that the point $z(c)$ is periodic for the parameter $c$
on $\gamma$. This analytic continuation is uniquely possible because 
we never encounter multipliers $+1$ for a period-$n$-orbit along
$\gamma$. Therefore, we obtain a unique periodic point $z(\tilde c)$
which is repelling. Since $\tilde c$ is within the connectedness
locus, the point $z(\tilde c)$ is the landing point of finitely many 
dynamic rays at periodic angles. 

At the parameter $c_0$, all the periodic dynamic rays of periods up
to $n$ land at repelling periodic points of periods up to $n$. These
periodic points can all be continued analytically along $\gamma$, they
remain repelling and keep their dynamic rays because the curve avoided
parameter rays of periods up to $n$ which make up wake boundaries.
Therefore, at the parameter $\tilde c$, there is no dynamic ray of
period $n$ available to land at $z(\tilde c)$, so that the rays landing
at this point must have periods $sn$, for some $s\geq 2$. Therefore,
at least $s$ rays must land at $z(\tilde c)$ (in fact, the number of
rays must be exactly $s$). Now let $\tilde U$ be the largest open
neighborhood of $\tilde c$ in which this periodic orbit can be continued
analytically as a repelling periodic point, retaining all its $s$
dynamic rays. This neighborhood is the {\em wake of the orbit} and is at
the heart of Milnor's discussion in \cite{MiOrbits}. This wake is
bounded by two parameter rays at periodic angles, landing together at a
parabolic parameter of ray period $sn$ and orbit period at most $n$.
Denote this landing point by $\tilde c_0$. Obviously, $\tilde U$ cannot
contain the hyperbolic component $W$, so $\tilde U$ must be contained
within the wake of $W$. 

The point $\tilde c_0$ is on the boundary of a hyperbolic component 
of period at most $n$ within the wake of $W$. The wake of this
component contains $\tilde U$ and thus $\tilde c$, so this component can
only be $W$ by the assumptions on $\tilde c$. It follows that $\tilde U$
is a sublimb of $W$ at a rational internal angle. 
\qed

\remark  
This result goes a long way towards proving local connectivity of the
Multibrot sets at boundary points of hyperbolic components, as was
pointed out to us by G.~Levin. 

\begin{corollary}[Trivial Fibers at Hyperbolic Component Boundaries]
\label{CorAlmostFiberTrivialHypComps} \lineclear
In every Multibrot set, the fiber of every boundary point of any
hyperbolic component is trivial, except possibly at the root of a
primitive component. The fiber of the root contains no point within
the limb of the component.
\end{corollary}
\proof
Any boundary point $c$ of a hyperbolic component can obviously be
separated from any point within the closure of the component and from
every sublimb which it is not a root of. This shows that fibers are
trivial for boundary points at irrational internal angles and at
co-roots (which exist only for $d\geq 3$), because neither have sublimbs
attached. If $c$ is the root of a hyperbolic component $W$, then it can
obviously be separated from any other point within the closure of $W$ or
from any rational sublimb, so its fiber contains no point within the
limb of the component. If $c$ is the point where the component $W$
bifurcates from $W_0$, then both arguments combine to show that the
fiber of $c$ is trivial.
\qed

The following corollary has been suggested by John Milnor.
It will be strengthened in Corollary~\ref{CorRootsDisconnect}.
\begin{corollary}[Roots Do Not Disconnect Limbs]
\label{CorRootsDontDisconnect} \lineclear 
The limb of any hyperbolic component is connected.
\qedd
\end{corollary}
\remark
Note that we define the wake of a hyperbolic component to be open, so
that the limb of the component (the intersection of the wake with the
Multibrot set) is a relatively open subset of $\M_d$. Sometimes, wake
and limb are defined so as to also contain the root. With that
definition, the corollary says that the root does not disconnect the
limb.
\proof
Let $W$ be the hyperbolic component defining the wake. Assume that its
limb consists of more than one connected component, and let $K$ be a
connected component not containing $W$. Any point within $K$ must then,
by Theorem~\ref{ThmNoIrratDecos}, be contained within some rational
sublimb of $W$, and then all of $K$ must be contained within the same
sublimb. But since the entire set $\M_d$ is connected and the limb is
obtained from $\M_d$ by cutting along a pair of parameter rays
landing at the root of the limb, it follows that the limb is
connected as well.
\qed

Here is another corollary which has been found independently by
Lavaurs~\cite[Proposition~1]{La}, Hubbard (unpublished) and
Levin \cite[Theorem~7.3]{Le}; another proof is in Lau and
Schleicher~\cite[Lemma~3]{IntAddr}. 

\begin{corollary}[Analytic Continuation Over Entire Wake]
\label{CorAnalytContWake}\lineclear  
If $c\in\M_d$ is a parameter which has an attracting periodic point
$z$, then this periodic point can be continued analytically as an
analytic map $z(c)$ over the entire wake of the hyperbolic component
containing $c$, and it is repelling away from the closure of the
component. Within any subwake at internal angle $p/q$, the point
$z(c)$ is the landing point of exactly $q$ dynamic rays with
combinatorial rotation number $p/q$. 
\end{corollary}
\remark
For parameters within a hyperbolic component of period $n$, there are
$d$ distinguished periodic points: the point on the attracting orbit
in the Fatou component containing the critical value, and $d-1$ further
boundary points of the same Fatou component with the same ray period
(and possibly smaller orbit period). If $n>1$, then exactly one of these
boundary points is the landing point of at least two dynamic rays (this
point is the ``dynamic root'' of the Fatou component), it is repelling
and can be continued analytically over the entire wake of the component,
retaining its dynamic rays (this is almost built into the definition of
the wake). The remaining $d-2$ distinguished boundary points of the
Fatou component are landing points of one dynamic ray each (they are
``dynamic co-roots''), and they also remain repelling and keep their
rays throughout the wake. For the attracting orbit, this is not true:
the rays landing at the orbit change, but the orbit nonetheless remains
repelling away from the closure of the hyperbolic component. 

\proof
Denote the hyperbolic component by $W$, let $c_0$ be its root and let
$n$ be the period. The orbit $z(c)$ can be continued analytically
over a neighborhood of $W-\{c_0\}$ within the wake of $W$, and it
will be repelling there (\cite[Section~7]{MiOrbits}, \cite[Section~5]
{ExtRays}). It will therefore remain repelling within any subwake of
$W$ at rational internal angles because that subwake is the wake where
our repelling orbit has certain dynamic rays at specific angles landing. 
Within the wake of $W$, away from the closure of the component and
outside of its rational subwakes, every Julia set is disconnected by
Theorem~\ref{ThmNoIrratDecos} and every orbit is repelling. Therefore,
$z(c)$ is repelling within the entire wake of $W$, except on the closure
of the component, and can be continued analytically throughout the
entire wake. The statement about the rays follows; compare for example
\cite{MiOrbits}. 
\qed
\remark
Traversing the wake of $W$ outside of $\M_d$, the combinatorial
rotation number of the orbit $z(c)$ behaves (locally) monotonically
with respect to external angles, rotating around $\Circle$ $d-1$
times. See also the appendices in \cite{MiOrbits} and \cite{GM}.

\newsection{Roots of Hyperbolic Components}
\label{SecRoots}

In this section, we will prove triviality of fibers also at roots of
primitive components. These roots are not handled by the arguments
from the previous section, or by the Yoccoz inequality. We start with
a combinatorial preparation from \cite{MandelStruct}.

\begin{lemma}[Approximation of Ray Pairs, Primitive Periodic Case]
\label{LemApproxPeriodic} \lineclear
Let $\theta<\theta'$ be the two periodic external angles of the
root of a primitive hyperbolic component. Then there exists a sequence
of rational parameter ray pairs $(\theta_n,\theta'_n)$ such that
$\theta_n\nearrow\theta$ and $\theta'_n\searrow\theta'$.
\end{lemma}
\sketch
There are various combinatorial or topological variants of this
proof. For example, it can be derived easily from the Orbit Separation
Lemma (see \cite[Lemma~\LemOrbitSep]{ExtRays} in the quadratic case
and \cite{Dominik} in the general case). We will use Hubbard trees.
They have been defined by Douady and Hubbard \cite{Orsay} for
postcritically finite polynomials as the unique minimal trees within
the filled-in Julia set connecting the critical orbit and traversing
any bounded Fatou components only along internal rays.

Consider the dynamics at the center of the primitive hyperbolic
component. The critical value is an endpoint of the Hubbard tree.
Since the component is primitive, there is a non-empty open subarc
$\gamma$ of the Hubbard tree avoiding periodic Fatou components which
ends at the Fatou component containing the critical value. We can
suppose this arc short enough so that it meets no branch point of the
Hubbard tree. Since the dynamics is expanding, the arc $\gamma$
contains a point $z$ which maps after some finite number $n$ of steps
onto the critical value. Choose $n$ minimal and restrict $\gamma$ to
the arc between $n$ and the periodic Fatou component containing the
critical value. Then the $n$-th iterate of the dynamics map $\gamma$
homeomorphically so that one end of $\gamma$ lands at the critical
value and the other end is a point on the critical orbit other than
the critical value. Therefore, $\gamma$ maps over itself in an
orientation reversing way and must contain some point $p$ which is
fixed under the $n$-th iterate. This point must be a repelling
periodic point, and it must be the landing point of two periodic
dynamic rays at external angles $\theta_n<\theta'_n$. This
construction can be done with arbitrarily short initial arcs $\gamma$,
producing periodic points arbitrarily close to the periodic Fatou
component containing the critical value. Alternatively, one can use
the fact that one endpoint of $\gamma$ is periodic (it is the landing
point of the dynamic rays at angles $\theta$ and $\theta'$), so that
it is possible to pull $\gamma$ back onto a subset of itself, creating
a preperiodic inverse image of $p$ closer to the critical value Fatou
component. Both ways, it is easy to see that the external angles of
the created periodic or preperiodic points have external angles
converging to $\theta$ and $\theta'$ from below respectively from
above.

Finally, Theorem~\ref{ThmRayCorrespondence} transfers these dynamic
ray pairs into parameter space.
\qed

\hide{
This lemma helps to finish the proof that the fiber of every boundary
point of a hyperbolic component is trivial.
\bigskip
\noindent
\begin{minipage}{\textwidth}
\begin{theorem}[Triviality of Fibers at Hyperbolic Components]
\label{ThmHypCompFiberTrivial} \lineclear
The fiber of every point on the closure of any hyperbolic component
of a Multibrot set $\M_d$ is trivial. 
\end{theorem}
\end{minipage}
\medskip
\proof
In Corollary~\ref{CorAlmostFiberTrivialHypComps}, we have already done
most of the work. It only remains to show that, if $c_0$ is the root of
a primitive hyperbolic component, then the fiber of $c_0$ contains no
point outside of the wake of the component.
Fix some primitive hyperbolic component and let $c_0$ be its root.
}

\begin{lemma}[Fiber of Primitive Root]
\label{LemFiberPrimRoot} \lineclear
Let $c_0$ be the root of any primitive hyperbolic component and let
$\tilde c\neq c_0$ be some point in $\M_d$ which is not in the limb
of $c_0$. Then there is a parameter ray pair at rational angles
separating $c_0$ from $\tilde c$ and from the origin. In
particular, the fiber of the root of a primitive hyperbolic
component contains no point outside of the wake of the component.
\end{lemma}
\proof
The strategy of the proof is simple: if we denote the periodic
parameter ray pair landing at $c_0$ by $(\theta,\theta')$, then 
Lemma~\ref{LemApproxPeriodic} supplies a sequence of parameter ray
pairs $(\theta_n,\theta'_n)$ converging to $(\theta,\theta')$,
and in the dynamics of $c_0$ there are characteristic dynamic ray
pairs at the same angles. \hide{Now let $\tilde c\neq c_0$ be a
parameter within the fiber of $c_0$.} If not all of these
characteristic dynamic ray pairs exist for $\tilde c$, then $\tilde
c$ is separated from $c_0$ by one of the parameter ray pairs
$(\theta_n,\theta'_n)$ (because these ray pairs bound the regions
in parameter space for which the corresponding dynamic ray pairs
exist). In this case, the two points $c_0$ and $\tilde c$ can be
separated as claimed and are thus in different fibers.

We may hence assume that all these dynamic ray pairs exist for the
parameter $\tilde c$. We will then show that the limiting dynamic
rays at angles $\theta$ and $\theta'$ also form a ray pair: this
forces $\tilde c$ to be in the closure of the wake of $c_0$,
contradicting our assumption.
\hide{, and we know by
Corollary~\ref{CorAlmostFiberTrivialHypComps} that the fiber of
$c_0$ does not contain points within this wake, while the only
boundary point of the wake within the Multibrot set is $c_0$. }
This finishes the proof of the lemma.

The only thing we need to prove, then, is the following claim:
{\em assume that for the parameter $\tilde c$, the dynamic rays at
angles $\theta_n$ and $\theta'_n$ land together for every $n$. Then the
dynamic rays at angles $\theta$ and $\theta'$ also land together.}

To prove this claim, denote the period of $\theta$ and $\theta'$ by
$k$. Since hyperbolic components and Misiurewicz points outside of the
wake of $c_0$ are landing points of rational parameter rays, they can
easily be separated from $c_0$ by the approximating parameter ray pairs,
and we may thus suppose that $\tilde c$ is neither a Misiurewicz point
nor on the closure of a hyperbolic component.

We modify the sequence $(\theta_n,\theta'_n)$ as follows: let
$(\theta_0,\theta'_0)$ be close enough to $(\theta,\theta')$ so that,
in the dynamics of $c_0$, these two ray pairs do not enclose a
postcritical point. Further, let the ray pairs $(\theta_n,\theta'_n)$ be
the unique ray pairs in between satisfying $\theta-\theta_n =
d^{-nk}(\theta-\theta_0)$ and $\theta'_n-\theta' =
d^{-nk}(\theta'_0-\theta')$. Denoting the regions bounded by
$(\theta,\theta')$ and $(\theta_n,\theta'_n)$ by $U_n$, then the
pull-back for $k$ steps of any $U_n$ along the periodic backwards orbit
of $(\theta,\theta')$ will map univalently onto $U_{n+1}$, so the same
branch of the pull-back fixing $(\theta,\theta')$ will map
$(\theta_n,\theta'_n)$ onto $(\theta_{n+1},\theta'_{n+1})$. 

All this works in the dynamics of $c_0$. For the parameter $\tilde
c$, denote the landing points of the two dynamic rays at angles
$\theta$ and $\theta'$ by $w$ and $w'$, respectively. Then the
critical value cannot be separated from $w$ or $w'$ by any rational
ray pair $(\alpha,\alpha')$: there is no such ray pair for $c_0$, and
the largest neighborhood of $\tilde c$ where such a ray pair exists
will be bounded by parameter rays at rational angles. Such a
neighborhood exists because $\tilde c$ is neither a Misiurewicz point
nor on the closure of a hyperbolic component. 

Similarly, for any $m>0$, the $m$-th forward image of the critical value
cannot be separated from the landing points of the dynamic rays at
angles $2^m\theta$ and $2^m\theta'$, or we could pull back such a
separating ray pair, obtaining a separation between the critical value
and the dynamic rays at angles $\theta$ and $\theta'$. It follows that
all the $(\theta_n, \theta'_n)$ (for $n\geq 1$) are separated from the
critical orbit by other such ray pairs with smaller or larger $n$. 

Now let $X$ be the complex plane minus the postcritical set for the
parameter $\tilde c$; then $X$ is the open set of all points in $\C$
which have a neighborhood that is never visited by the critical orbit.
Let $p(z)=z^d+\tilde c$ and set $X':=p^{-1}(X)$; since the postcritical
set is forward invariant, we have $X'\subset X$, and $p\colon X'\to X$
is an unbranched covering, i.e., a local hyperbolic isometry. 

Let $u,u'$ be two points on the dynamic rays at angles $\theta$ and
$\theta'$, respectively, and construct a simple piecewise analytic
curve $\gamma_0$ as follows: connect $u$ along an equipotential to the
dynamic ray at angle $\theta_0$ (the ``short way'', decreasing
angles), then continue along this ray towards its landing point, then
out along the dynamic ray at angle $\theta'_0$ up to the potential of
$u'$, and connect finally along this equipotential to $u'$
(increasing angles). This curve runs entirely within $X$, so it has
finite hyperbolic length $\ell_0$, say.

Now we pull this curve back along the periodic backwards orbit of
$\theta$ and $\theta'$, yielding a sequence of curves $\gamma_1,
\gamma_2,\ldots$ after $k, 2k,\ldots$ steps. The branch of the
pull-back fixing $\theta$ will also fix $\theta'$: we had verified
that in the dynamics of $c_0$, and this is no different for $\tilde
c$ because the critical value is always on the same side of all the
ray pairs.

Denote the hyperbolic lengths of $\gamma_n$ in $X$ by $\ell_n$.
Since $X'\subset X$, this sequence is strictly monotonically
decreasing. The two endpoints of the curves in the sequence
obviously converge to the two landing points $w$ and $w'$ of the
dynamic rays at angles $\theta$ and $\theta'$. Now there are two
possibilities: either the sequence $(\gamma_n)$ has a subsequence
which stays entirely within a compact subset of $X$, or it does
not. If it does, then the hyperbolic length shrinks by a definite
factor each time the curve is within the compact subset of $X$, so
that we can connect points (Euclideanly) arbitrarily closely to $w$
and $w'$ by (hyperbolically) arbitrarily short curves at bounded
(Euclidean) distances from the postcritical set, and this implies
$w=w'$ as required. On the other hand, if these curves converge to
the boundary, then their Euclidean lengths must shrink to zero, and
we obtain the same conclusion. The proof of the lemma is complete.
\qed

\begin{corollary}[Triviality of Fibers at Hyperbolic Components]
\label{CorHypCompFiberTrivial} \lineclear
The fiber of every point on the closure of any hyperbolic component
of a Multibrot set $\M_d$ is trivial. 
\end{corollary}
\proof
In Corollary~\ref{CorAlmostFiberTrivialHypComps}, we have already
done most of the work; the remaining statement is exactly the
content of the previous lemma.
\qed

The next corollary strengthens
Corollary~\ref{CorRootsDontDisconnect} and has also been suggested
by Milnor.
\begin{corollary}[Roots Disconnect Multibrot Sets]
\label{CorRootsDisconnect} \lineclear
Every root of a hyperbolic component of period greater than $1$
disconnects its Multibrot set into exactly two connected components.
\end{corollary}
\proof
Let $c_0$ be the root of a hyperbolic component of period greater
than $1$. It is the landing point of exactly two parameter rays at
periodic angles. This parameter ray pair disconnects $\C$ into
exactly two connected components. The component not containing the
origin is the {\em wake} of the component; its intersection with
$\M_d$ is the limb. The limb is connected by
Corollary~\ref{CorRootsDontDisconnect}. Now let $\tilde c\in\M_d$
with $\tilde c\neq c_0$ be any parameter not in the limb of $c_0$. 

We first discuss the case that the component is primitive. By
Lemma~\ref{LemFiberPrimRoot} above, there is a rational parameter
ray pair $S$ separating $c_0$ from $\tilde c$ and from the origin.
This ray pair $S$ disconnects $\C$ into two connected components.
Let $\M'$ be the closure of the connected component of the origin
intersected with $\M_d$. Then $\M'$ itself is connected and
contains the origin and $\tilde c$. It follows that $\tilde c$ and
$0$ are in the same connected component of $\M_d-\{c_0\}$, so $c_0$
cuts $\M_d$ into exactly two connected components.

In the non-primitive case, the parameter $c_0$ is the root of one
hyperbolic component and on the boundary of another one, say $W_0$.
Then by Theorem~\ref{ThmNoIrratDecos}, there are three cases: the
point $c_0$ may be outside of the limb of $W_0$ so that the two
parameter rays landing at the root of $W_0$ separate $c_0$ and
$\tilde c$; the point $c_0$ may be on the closure of $W_0$; or it
may be in a sublimb of $W_0$ at rational internal angles, but not
in the sublimb with $c_0$ on its boundary (which is the wake of
$c_0$). In all three cases, it is easy to see that $\tilde c$ must
be in the same connected component of $\M_d-\{c_0\}$ as the origin. 

Therefore, $c_0$ disconnects $\M_d$ into exactly two connected
components in the primitive as well as in the non-primitive case.
\qed

\begin{corollary}[Rays at Boundary of Hyperbolic Components]
\label{CorRaysAtHypComps} \lineclear
Every boundary point of a hyperbolic component at irrational internal
angle is the landing point of exactly one parameter ray, and the
external angle of this ray is irrational (and in fact transcendental). 

Every boundary point at rational internal angle is the landing point
of exactly two parameter rays, except co-roots: these are the landing
points of exactly one parameter ray. The external angles of all these
rays are periodic and in particular rational. 

In no case is a boundary point of a hyperbolic component in the
impression of any further parameter ray.
\end{corollary}
\proof
Every boundary point of a hyperbolic component at irrational internal
angle is in the boundary of $\M_d$ and thus in the impression of some
ray. Since the fiber of this boundary point is trivial, the ray must
land there and its impression is a single point
\cite[Lemma~\LemImpressFiber]{Fibers}. And since all the parameter
rays at rational rays land elsewhere, the external angle of the ray
must be irrational (and in fact transcendental; see \cite{BS} for a
proof in the quadratic case). If this point is in the impression of a
further ray, this ray must land there, too, and these two rays
separate some open subset of $\C$ from the component. Since not all
rays between the two landing rays can land at the same point (the
rational rays land elsewhere), there must be some part of $\M_d$
between these two rays, and this contradicts
Theorem~\ref{ThmNoIrratDecos}. 

We know that boundary points at rational internal angles are
parabolic points, and the statements about the landing properties of
rational parameter rays are well known (compare
Section~\ref{SecMultibrot}). If such a point is in the impression of
an irrational parameter ray and thus its landing point, we get a
similar contradiction as above, except if the parabolic point is the
root of a primitive component and the extra parameter ray is outside
of the wake of the component. But that case is handled conveniently by
Lemma~\ref{LemApproxPeriodic}, supplying lots of parameter ray pairs
which separate any parameter ray outside of the wake of the component
from the component, its root, and all of its co-roots. 
\qed

\newsection{Misiurewicz Points}
\label{SecMisiurewicz}

Now we turn to Misiurewicz points, proving triviality of fibers in a
somewhat similar way. Again, we need some combinatorial preparations
from \cite{MandelStruct}. They are preperiodic analogs to
Lemma~\ref{LemApproxPeriodic}. 

\begin{lemma}[Preperiodic Critical Orbits Have Trivial Fibers]
\label{LemMisiuDynamics} \lineclear
Let $c_0$ be a Misiurewicz point and let $z$ be the any periodic point
on the forward orbit of the critical value. Then there are finitely
many rational ray pairs separating $z$ from the entire critical orbit
except $\{z\}$. The separating ray pairs can be chosen so that their
landing point is not on the grand orbit of the critical point.
\end{lemma}

\sketch
Like Lemma~\ref{LemApproxPeriodic}, this result has various proofs.
For example, it is a consequence of
\cite[Theorem~\ThmPeriodicFiber]{Fibers}. Another variant uses the
fact that the Julia set is locally connected, so its fibers are
trivial when constructed using rational ray pairs 
\cite[Proposition~\PropLocConnJuliaFiber]{Fibers}.

Our proof uses similar ideas as in Lemma~\ref{LemApproxPeriodic}.
Again, we consider the Hubbard tree of the Julia set (the unique
minimal tree connecting the critical orbit within the Julia set). 
It suffices to prove the result for any point on the periodic orbit of
$z$. 

\hide{
, so we may suppose $z$ is the {\em characteristic} point on this
orbit: this is the unique periodic point the critical value can be
connected to within the Julia set without crossing any other point on
the same periodic orbit (the existence of this characteristic point is
at the heart of many combinatorial descriptions of quadratic
polynomials: see e.g.\ \cite{Thurston}, \cite{MiOrbits} or
\cite{ExtRays}; the proofs also work for arbitrary degree unicritical
polynomials \cite{Dominik}).
}

Let $\gamma$ be a short arc in the Hubbard tree starting at $z$
towards the critical value $c$. By expansion, there must be a point
$z_n$ on $\gamma$ which maps to $c$ after some finite number $n$ of
iterations, and we may suppose $n$ minimal. Restrict $\gamma$ to the
arc between $z$ and $z_n$. The $n$-th iterate will then map $\gamma$
forward homeomorphically such that a subarc of $\gamma$ maps over
$\gamma$. If $n$ is not a multiple of the period of $z$, then the
$n$-th iterate does not fix any end of $\gamma$, and there must be an
interior point $p$ of $\gamma$ which is fixed under this iterate. If
$n$ is a multiple of the period of $z$, the $n$-th iterate of $\gamma$
connects $z$ to $c$ homeomorphically. We keep mapping this arc forward
until is crosses the critical point for the first time. The next
iterate will then provide a periodic point $p$ unless the point $z$ is
back again. If this happens all the time, then the itineraries of $z$
and $c$ must coincide, but this is not the case: the first is
periodic, the second preperiodic. 

Once we have periodic points close to $z$, we have them at all
branches of $z$ because local branches at periodic points are permuted
transitively, except possibly if there are exactly two branches which
are fixed by the first return map of $z$ (compare
\cite[Lemma~\MiLem]{MiOrbits} or \cite[Lemma~\ExtRayLem]{ExtRays}).
In that case, we have periodic points close to any point on the
periodic orbit of $z$ at least from the side towards the critical
point, and it is not hard to verify that the dynamics cannot always
preserve these sides. Therefore, mapping forward along the periodic
orbit, we get periodic points arbitrarily close to $z$ from all sides.
If the approximating periodic points are chosen close enough to $z$,
they cannot be from the periodic orbit of $z$.
\qed

\begin{lemma}[Approximation of Ray Pairs, Preperiodic Case]
\label{LemApproxPreperiodic} \lineclear
Consider a Misiurewicz point and let\/ $\{\theta_i\}$ be the finite
set of its preperiodic external angles (i.e., the external angles of
the critical value in its dynamic plane). Then, for every
$\eps>0$, there are finitely many rational parameter ray pairs
separating the Misiurewicz point from every parameter ray whose
external angle differs at least by $\eps$ from all of the $\theta_i$:
one can take one such rational ray pair for every connected component
of\/ $\Circle-\{\theta_i\}$.  All these ray pairs can be chosen to be
periodic.
\end{lemma}
\sketch
From the previous lemma, we can conclude that any periodic point $z$
on the critical orbit can be separated from the remaining critical
orbit by periodic dynamic ray pairs. These ray pairs can be pulled
back along the repelling periodic orbit of $z$, at least if they were
close enough to $z$, and we get dynamic ray pairs arbitrarily close
to the rays landing at the critical value. We need to transfer them
into parameter space. This is easily done using
Theorem~\ref{ThmRayCorrespondence} for the rays separating the
critical value from the critical point. For the others, it requires a
more precise description of the combinatorics of Multibrot sets; 
see \cite{MandelStruct}. The idea is to use 
Theorem~\ref{ThmRayCorrespondence} not in the dynamic plane of the
Misiurewicz point, but of hyperbolic components in the various wakes
of the Misiurewicz point. A sketch of the quadratic case is in
Douady~\cite{DoCompacts}.  
\qed

\remark
The statement essentially says that every external ray which is not
landing at a given Misiurewicz point can be separated from this point
by a rational ray pair. The approximating ray pairs can be required
at will to land at Misiurewicz points, or at roots of primitive
hyperbolic components, or at roots of non-primitive hyperbolic
components (in the latter case, such a non-primitive hyperbolic
component will usually be a bifurcation by a factor two). See
\cite{MandelStruct}.

The following obvious corollary is an analogue to
Corollary~\ref{CorRaysAtHypComps}.
\begin{corollary}[No Irrational Rays at Misiurewicz Points]
\label{CorNoIrratMisiuRays} \lineclear
No Misiurewicz point is in the impression of a parameter ray at an
irrational external angle.
\qedd
\end{corollary}

The following theorem is well known for the Mandelbrot set. As far as
I know, it was first proved by Yoccoz using his puzzle techniques. A
variant of this proof for most cases is indicated in Hubbard's
paper \cite[Theorem~14.2]{HY}. Tan Lei has published another proof in
\cite{TLLocConn}.

\begin{theorem}[Misiurewicz Points Have Trivial Fibers]
\label{ThmMisiuFiber} \lineclear
The fiber of any Misiurewicz point in any Multibrot set is trivial.
\end{theorem}
\proof
Let $c_0$ be a Misiurewicz point and denote its various subwakes by
$U_0,U_1,\ldots,U_s$, for some integer $s\geq 0$, such that $U_0$
contains the origin. Let $\theta_0,\ldots,\theta_s$ be the external
angles of the Misiurewicz point. We have to show that, for every
$\tilde c\in U_i\cap\M_d$, there is a rational ray pair separating
$\tilde c$ from $c_0$. 

If $\tilde c$ is on the closure of a hyperbolic component, then one
of the approximating ray pairs will separate the entire hyperbolic
component from $c_0$, and we are done. Similarly, we are done if
$\tilde c$ is a Misiurewicz point. If $\tilde c$ is in some
non-hyperbolic component (if that ever occurs), then we may replace
$\tilde c$ by some other point (different from $c_0$) on the boundary
of the same non-hyperbolic component and thus on the boundary of the
Multibrot set: a ray pair separating a boundary point of a
non-hyperbolic component from $c_0$ will separate the entire component
from $c_0$. We may therefore suppose that $\tilde c$ is a boundary point
of the Multibrot set, that all periodic orbits are repelling, and
that the critical orbit is infinite.

In the dynamic plane of $c_0$, let $z$ be the first periodic point
on the orbit of the critical value. In a parameter neighborhood of
$c_0$, this periodic point can be continued analytically as a
periodic point $z(c)$. We know from Lemma~\ref{LemMisiuDynamics} that
there are finitely many rational ray pairs separating $z(c_0)$ from
its entire forward orbit, such that the neighborhood of $z(c_0)$ cut
out by these rays together with any equipotential, when pulled back
along the periodic backwards orbit of $z(c_0)$, shrinks to the point
$z(c_0)$ alone. The same will then be true for the point $z(c)$ for
parameters $c$ from a sufficiently small neighborhood $V$ of $c_0$ 
because the finitely many bounding ray pairs will persist under
perturbations, and the critical value will not be contained in the
regions to be pulled back. The point $z(c)$ can thus be separated from
any other point in the Julia set by a rational ray pair (in other
words, the fiber of $z(c)$ is hence trivial when fibers of the Julia
set are constructed using rational external rays. Another way to
arrive at this same conclusion is to use 
\cite[Theorem~\ThmPeriodicFiber]{Fibers}.)

Let $z'(c)$ be the analytic continuation of the repelling preperiodic
point which equals the critical value for the parameter $c_0$. This
analytic continuation is possible sufficiently closely to $c_0$;
we may assume it to be possible throughout $V$, possibly by shrinking
$V$. Since the fiber of $c_0$ is connected, we may suppose that $\tilde
c\in V$, subject to the same restrictions as above, so that $z(\tilde
c)$ and thus $z'(\tilde c)$ have trivial fibers. 

Since $z'(\tilde c)$ is different from the critical value $\tilde c$,
there is a rational ray pair $(\tilde \theta,\tilde \theta')$
separating the critical value from $z'(\tilde c)$ and from all the
dynamic rays at angles $\theta_i$, which will land at $z'(\tilde c)$.
The angles of this separating ray pair will differ by some $\eps>0$ from
all the $\theta_i$. 

This separating dynamic ray pair will persist under sufficiently
small perturbations of $\tilde c$. Since we had assumed $\tilde c \in
\partial\M_d$, such a perturbation is possible into the exterior of
$\M_d$. The external angle of the perturbed parameter must then
differ by at least $\eps$ from all the $\theta_i$, no matter how
small the perturbation. Now we use Lemma~\ref{LemApproxPreperiodic}
to get finitely many rational parameter ray pairs separating $c_0$
from all the parameter rays the angles of which differ from all the
$\theta_i$ by at least $\eps/2$. This separation must then also
separate $\tilde c$ from $c_0$, so these two parameters are in different
fibers. It follows that the fiber of $c_0$ is trivial.
\qed

\begin{corollary}[Misiurewicz Points Disconnect]
\label{CorMisiuDisconnect} \lineclear
Every Misiurewicz point disconnects its Multibrot set in exactly as
many connected components as there are rational parameter rays landing
at it.
\end{corollary}
\proof
All the parameter rays landing at a Misiurewicz point have preperiodic
external angles by Corollary~\ref{CorNoIrratMisiuRays}, and they
obviously disconnect the Multibrot set in at least as many connected
components as there are such rays.

By Lemma~\ref{LemApproxPreperiodic}, between any two adjacent
external rays of a Misiurewicz point there is a collection of
rational parameter ray pairs exhausting the interval of external
angles in between, so any extra connected component at the
Misiurewicz point cannot have any external angles. It must therefore
be contained within the fiber of the Misiurewicz point, but this
fiber is the Misiurewicz point alone.
\qed

\newsection{Fibers and Combinatorics}
\label{SecCombinatorics}

Now that we have trivial fibers at all the landing points of rational
rays, it follows that fibers of any two points of $\M_d$ are either
equal or disjoint. This is Lemma~{\LemFibersNice} in \cite{Fibers}; we
repeat the proof in order to make this paper self-contained.

\begin{theorem}[Fibers of $\M_d$ are Equivalence Classes]
\label{ThmFibersNice} \lineclear
The fibers of any two points in $\M_d$ are either equal or disjoint.
\end{theorem}
\proof
The relation ``$c_1$ is in the fiber of $c_2$'' is always reflexive,
and it is symmetric for $\M_d$ by Lemma~\ref{LemFibers}. In order to
show transitivity, assume that two points $c_1$ and $c_2$ are both in
the fiber of $c_0$. If they are not in the fibers of each other, then
the two points can be separated by a separation line avoiding $c_1$ and
$c_2$. If such a separation line can avoid $c_0$, then these two
points cannot both be in the fiber of $c_0$. The only separation
between $c_1$ and $c_2$ therefore runs through the point $c_0$, so
$c_0$ cannot be in the interior of $\M_d$ and rational rays land at
$c_0$. Therefore, the fiber of $c_0$ consists of $c_0$ alone. Any
two points with intersecting fibers thus have indeed equal fibers. 
\qed

The theorem allows to simply speak of {\em fibers} of $\M_d$ as
equivalence classes of points with coinciding fibers, as opposed to
``fibers of $c$'' for $c\in\M_d$. Another consequence shown in 
Lemma~{\LemFibersNice} in \cite{Fibers} is that there is an obvious
map from external angles to fibers of $K$ via impressions of external
rays. This map is surjective onto the set of fibers meeting $\partial
K$.

The following corollary is obvious and just stated for easier
reference. In Corollary~\ref{CorMandelInterior}, we had only been able
to show that at most one parameter ray at a rational angle can land at
a non-hyperbolic component.

\begin{corollary}[Non-Hyperbolic Components Rationally Invisible]
\label{CorNoNonHypRatRay} \lineclear
No rational parameter ray lands on the boundary of a non-hyperbolic
component.
\qedd
\end{corollary}

Combinatorial building blocks of the Multibrot sets which are often
discussed are combinatorial classes. We will show that they are
closely related to fibers.

\begin{definition}[Combinatorial Classes and Equivalence]
\label{DefCombClass} \lineclear 
We say that two connected Julia sets are\/ {\em combinatorially
equivalent} if in both dynamic planes external rays at the same
rational angles land at common points. Equivalence classes under this
relation are called\/ {\em combinatorial classes}.
\end{definition}

\remark 
In the language of Thurston~\cite{Th}, combinatorially equivalent Julia
sets are those having the same rational lamination. The definition is
such that topologically conjugate (monic) Julia sets are also
combinatorially equivalent (for an appropriate choice of one of the
$d-1$ fixed rays to have external angle $0$). In particular, all the
Julia sets within any hyperbolic component, at its root and at its
co-roots are combinatorially equivalent.  

With the given definition, the combinatorial class of a hyperbolic
component also includes its boundary points at irrational angles,
although the dynamics will be drastically different there. Therefore,
one might want to refine the definition of combinatorial equivalence
accordingly. 

\remark
\looseness 1
There are homeomorphic Julia sets which are not topologically
conjugate and which are thus not combinatorially equivalent: as an
example, it is not hard to see that the Julia set of $z^2-1$ (known
as the ``Basilica'') is homeomorphic to any locally connected
quadratic Julia set with a Siegel disk of period one. However, this
homeomorphism is obviously not compatible with the dynamics, and the
Basilica is in a different combinatorial class than any Siegel disk
Julia set.

\begin{proposition}[Combinatorial Classes and Fibers]
\label{PropClassFiber} \lineclear 
Hyperbolic components together with their roots, co-roots and
irrational boundary points form combinatorial classes. All other
combinatorial classes are exactly fibers. In particular, if there is
any non-hyperbolic component, then its closure is contained in a
single combinatorial class and a single fiber.
\end{proposition}
\proof 
\looseness 1
The landing pattern of dynamic rays changes upon entering the wake of
a hyperbolic component, so any parameter ray pair at periodic angles
separates combinatorial classes. Similarly, the landing pattern is
different within all the subwakes of any Misiurewicz point: in a
neighborhood of the Misiurewicz point, the rays landing at the critical
value will have the same angles, but the $d$ preperiodic inverse images
will carry different angles according to which subwake of the
Misiurewicz point the parameter is in. Hence parameter ray pairs at
preperiodic angles also separate combinatorial classes, and
combinatorial classes are just what fibers would be if separations were
allowed only by rational ray pairs, excluding separation lines
containing curves through interior components. Since fibers have more
separation lines, they are contained in combinatorial classes. The
closure of any non-hyperbolic component is contained in a single fiber
(Corollary~\ref{CorMandelInterior}), so it is also contained within a
single combinatorial class. 

We know that the rational landing pattern is constant throughout
hyperbolic components and at its root and co-roots, as well as at its
irrational boundary points. Hyperbolic components with these specified
boundary points are therefore in single combinatorial classes; since
every further point in the Multibrot set is either in a rational subwake
of the component or outside the wake of the component, the
combinatorial class of any hyperbolic component is exactly as
described. 

Finally, we want to show that every non-hyperbolic combinatorial
class is a single fiber. If this was not so, then two points within
the combinatorial class could be separated by a separation line.
Unless these points are both on the closure of the same hyperbolic
component, such a separation line can always be chosen as a ray pair
at rational angles, and we have seen above that such ray pairs bound
combinatorial classes.
\qed

\newpage
\typeout{newpage}

The following corollary is closely related to the Branch
Theorem~\ref{ThmBranch}.

\begin{corollary}[Three Rays at One Fiber]
\label{CorThreeRaysOneFiber} \lineclear 
If three parameter rays of a Multibrot set accumulate at the same
fiber, then the three rays are preperiodic and land at a common
Misiurewicz point. If three parameter rays accumulate at a common
combinatorial class, then this combinatorial class is either a
Misiurewicz point or a hyperbolic component, and the three rays
land. 
\end{corollary}
\proof
We will argue similarly as in Corollary~\ref{CorLocConnFiberHyp}.
Denote the three external angles by $\theta_1,\theta_2,\theta_3$ in
increasing order and let $Y$ be their common fiber. These angles
separate $\Circle$ into three open intervals. Two of these intervals
do not contain the angle $0$. Choose a preperiodic rational angle
from both of them. The corresponding parameter rays land at two
Misiurewicz points. Applying the Branch Theorem~\ref{ThmBranch}
to these two points, we find either that one of these two points
separates the other from the origin, or there is a Misiurewicz point
or hyperbolic component which separates both from each other and from
the origin. In all cases, it is impossible to connect the three
parameter rays at angles $\theta_{1,2,3}$ to a single fiber, unless
these three rays land at the separating Misiurewicz point or
hyperbolic component. In the Misiurewicz case, all three angles must
be preperiodic by Corollary~\ref{CorNoIrratMisiuRays}, and in the
hyperbolic case the three parameter rays must still land at a common
point because of fiber triviality on closures of hyperbolic components.
This is impossible. 

If, however, the three rays are only required to land at a common
combinatorial class, then this combinatorial class may either be a
Misiurewicz point or a hyperbolic component, and both cases obviously
occur.
\qed

We have defined fibers of a Multibrot set $\M_d$ using parameter rays
at periodic and preperiodic angles (where the dynamics as usual
multiplication by the degree $d$). However, it turns out that only
periodic angles are necessary.

\begin{proposition}[Fibers Using Periodic Parameter Rays]
\label{PropPeriodicParaRays} \lineclear
Fibers for a Multibrot set $\M_d$ remain unchanged when they are
defined using only parameter rays at periodic angles, rather than at
all rational angles.
\end{proposition}
\proof
Preperiodic parameter rays never land on the boundary of an interior
component of a Multibrot set, so separation lines using preperiodic
rays are always ray pairs. We only have to show that if any two
parameters can be separated by a preperiodic ray pair, then there is a
periodic ray pair separating them as well. 

First we show that the fiber of any Misiurewicz point is trivial even
when only periodic parameter rays are used. Indeed, let $c_0$ be a
Misiurewicz point and let $c\in\M_d-\{c_0\}$ be a different parameter.
By triviality of the fiber of $c_0$, there is a separation line between
$c_0$ and $c$. If this separation line is a preperiodic ray pair, then
there is a periodic ray pair separating $c_0$ from the preperiodic ray
pair (see Lemma~\ref{LemApproxPreperiodic}) and thus also from $c$.
Therefore, $c_0$ can be separated from any parameter in $\M_d-\{c_0\}$
even if only parameter rays at periodic angles are allowed in the
construction of fibers.

Finally, let $c_1$ and $c_2$ be two parameters of $\M_d$ in different
fibers. If they are separated by a parameter ray pair at preperiodic
angles, let $c_0$ be the Misiurewicz point at which this ray pair
lands. Then there is a periodic ray pair separating $c_0$ and $c_1$.
Since the rays landing at $c_0$ separate $c_1$ and $c_2$, this
periodic ray pair must also separate $c_1$ and $c_2$ and we are done.
\qed

It follows that combinatorial classes are exactly the pieces that we
can split Multibrot sets into when using only periodic ray pairs as
separation lines (except that we have to declare what we want
combinatorial classes to be on the boundary of hyperbolic components).
This is exactly the partition used to define internal addresses, which
was introduced in \cite{IntAddr} mostly for postcritically finite
parameters. Now we see that combinatorial classes are the objects
where internal addresses live naturally. The following is a
restatement of \cite[Theorem~9.2]{IntAddr} and the remark thereafter. 

\begin{corollary}[Internal Addresses Label Combinatorial Classes]
\label{CorIntAddrCombClass} \lineclear
Angled internal addresses label combinatorial classes completely:
points within different combinatorial classes have different angled
internal addresses, and points within the same combinatorial class
have the same angled internal address.
\qedd
\end{corollary}

In order to label fibers, we only have to distinguish points within
hyperbolic components. This is most naturally done by the multiplier
of the unique non-repelling cycle. For degrees $d>2$, every
hyperbolic component contains $d-1$ points with equal non-repelling
cycles. They are distinguished by their ``sectors'' within the
component, and these sectors are reflected in both the kneading
sequences and the angled internal addresses with sectors 
\cite[Section~10]{IntAddr}.

\begin{theorem}[The Pinched Disk Model of Multibrot Sets]
\label{ThmPinchedDisk} \lineclear
The quotient of any Multibrot set in which every fiber is collapsed to
a point is a compact connected locally connected metric Hausdorff
space. In particular, the quotient is pathwise connected.
\end{theorem}
\proof
Every fiber of $\M_d$ is closed. In fact, the entire equivalence
relation is closed: suppose that $(z_n)$ and $(z'_n)$ are two
converging sequences in $\M_d$ such that $z_n$ and $z'_n$ are in a
common fiber for every $n$. The limit points must then also be in a
common fiber: if they are not, then they can be separated by either a
ray pair at periodic angles or by a separation line through a
hyperbolic component, and the separation runs in both cases only
through points with trivial fibers. In order to converge to limit
points on different sides of this separation line, all but finitely
many points of the two sequences must be on the respective sides of
the separation line, and $z_n$ and $z'_n$ cannot be in the same fiber
for large $n$. 

It follows that the quotient space is a Hausdorff space with respect
to the quotient topology. It is obviously compact and connected, and
local connectivity is easy: the same proof applies which shows that
trivial fibers imply local connectivity (see
Proposition~\ref{PropFiberLocConn}): it makes sense to speak of fibers
of the quotient space, and they are all trivial.

The quotient inherits a natural metric from $\C$: the distance between
any two points is the distance between the corresponding fibers. Since
fibers are closed, different points always have positive distances.
Symmetry and the triangle inequality are inherited from $\C$.
As a compact connected locally connected metric space, the quotient
space is pathwise connected (see Milnor~\cite[\S16]{MiIntro}).
\qed

This quotient space is called the ``pinched disk model'' of the
Multibrot set; compare Douady~\cite{DoCompacts}. It comes with an
obvious continuous projection map $\pi$ from the actual Multibrot set,
and inverse images of points under $\pi$ are exactly fibers of $\M_d$. 
The name ``pinched disk model'' comes from Thurston's lamination model
\cite{Th}: start with the closed unit disk $\diskbar$; for any
periodic ray pair at angles $\theta$ and $\theta'$, connect the
corresponding boundary points of $\diskbar$, for example by a
geodesic with respect to the hyperbolic metric in $\disk$ (or simply
by a Euclidean straight line). This connecting line is called a
``leaf''. Whenever two boundary points $\theta,\theta'$ are limits of
two sequences $(\theta_n)$, $(\theta'_n)$ such that
$(\theta_n,\theta'_n)$ is a leaf for every $n$, then the boundary
points $\theta$ and $\theta'$ in $\diskbar$ have to be connected as
well to form a leaf. Now we take the quotient of $\diskbar$ in which
each leaf is collapsed to a point: the disk $\diskbar$ is ``pinched''
along this leaf. The quotient space is homeomorphic to our space
constructed above.

\bigskip
\noindent Dierk Schleicher\\
\medskip
%\parbox[t]{65mm}{ 
Fakult\"at f\"ur Mathematik\\ Technische
Universit\"at\\ Barer Stra{\ss}e 23\\ D-80290 M\"unchen, Germany\\
{\sl dierk$@$mathematik.tu-muenchen.de} 
%}\hfill and \hfill
%\parbox[t]{65mm}{ Institute for Mathematical Sciences \\ State
%University of New York \\ Stony Brook, NY 11794-3660\\ {\sl
%dierk$@$math.sunysb.edu} }

%\vfill
%\hfill Revision of: 31.~Oktober 1998

\end{document}